\documentclass[pdf, ssy, preprint]{imsart}
\RequirePackage[OT1]{fontenc}
\RequirePackage{amsthm,amsmath}
\RequirePackage[colorlinks,citecolor=blue,urlcolor=blue]{hyperref}

\usepackage{float}
\usepackage{amsmath}
\usepackage{graphicx}
\usepackage{latexsym}
\usepackage{amsfonts}
\usepackage{amssymb}
\usepackage{amsthm}
\usepackage{amsfonts}
\theoremstyle{plain}
\newtheorem{theorem}{Theorem}
\newtheorem{lemma}{Lemma}
\newtheorem{remark}{Remark}
\newtheorem{cor}{Corollary}
\theoremstyle{definition}

\newcommand{\I}{{\bf I}}

\begin{document}
\begin{frontmatter}
\title{Stability of a Markov-modulated Markov Chain, with application to a wireless
network governed by two protocols}
\runtitle{Multi-radio}
\begin{aug}
\author{\fnms{Sergey} \snm{Foss} \thanksref{t1,m1,m2} \ead[label=e1]{s.foss@hw.ac.uk}},
\author{\fnms{Seva} \snm{Shneer} \thanksref{t1,m1} \ead[label=e2]{v.shneer@hw.ac.uk}}
\and
\author{\fnms{Andrey} \snm{Tyurlikov} \thanksref{t1,m3} \ead[label=e3]{turlikov@vu.spb.ru}}

\thankstext{t1}{The research of the authors
was partially supported by the Royal Society International Joint Project}

\affiliation{Heriot-Watt University, Edinburgh\thanksmark{m1}, Institute of Mathematics, Novosibirsk\thanksmark{m2} and St. Petersburg State University of Aerospace Instrumentation\thanksmark{m3}}

\address{School of Mathematical and Computer Sciences,\\ Heriot-Watt University, \\ EH14 4AS, Edinburgh, UK\\
\printead{e1}\\
\phantom{E-mail:\ }\printead*{e2}}
\address{St. Petersburg \\ State University of Aerospace Instrumentation\\
\printead{e3}}
\end{aug}

\date{\today}

\begin{abstract}
We consider a discrete-time Markov chain $(X^t,Y^t)$, $t=0,1,2,\ldots$,
where the $X$-component forms a Markov chain itself. Assume that
$(X^t)$ is Harris-ergodic and consider an auxiliary Markov chain $\{ \widehat{Y}^t\}$ whose transition probabilities
are the averages of transition probabilities of the $Y$-component of the $(X,Y)$-chain, where
the averaging is weighted by the stationary distribution of the $X$-component. 

We first provide natural conditions in terms of test functions ensuring that the $\widehat{Y}$-chain is positive recurrent and then prove that these conditions are also sufficient for positive recurrence of the original chain $(X^t,Y^t)$. The we prove a ``multi-dimensional''
extension of the result obtained. 
In the second part of the paper, we apply our results to two versions of a multi-access
wireless model governed by two randomised protocols.
\end{abstract}



\end{frontmatter}
\section{Introduction}

We develop an approach to the stability analysis based
on an averaging Lyapunov criterion.
More precisely, we consider a discrete-time Markov chain $\{Z^t, t=0,1,2,\ldots \}$
with values in a general state space which has two components, $Z^t=(X^t,Y^t)$.
We assume that the first component
 $\{ X^t\}$ forms a
Markov chain itself. We further assume that the
Markov chain $\{X^t\}$  is Harris-ergodic, i.e. there exists a unique
stationary distribution $\pi_X$ and, for any initial value $X^0=x$,
the distribution of $X^t$ converges to the stationary distribution
in the total variation norm:
$$
\sup_{A\in \mathcal {B_X}} |{\mathbf P} (X^t \in A) - \pi_X(A) |
\to 0 \quad \mbox{as} \quad t\to\infty.
$$

In many known models in applied probability, the time evolution of the process describing the state of the system may be expressed as a multi-component Markov chain where one of the components is either a Markov chains itself (see, e.g. \cite{Gamarnik}, \cite{Gamarnik2}, \cite{Litvak}, \cite{Borst} and references therein) or, more generally, behaves as``almost stable", changing only relatively slowly, 
while the other component changes more rapidly (see, e.g., \cite{Shah}).

Our aim is to formulate and prove a stability criterion for the
two-component Markov chain $Z^t= (X^t,Y^t)$ by making links to
an auxiliary Markov chain $(\widehat{X}^t,\widehat{Y}^t)$, where
$\{\widehat{X}^t\}$ is an i.i.d. sequence with common distribution
$\pi_X$.

First, one can recall a standard approach to the stability analysis
which may be used for this model and which is
based on the following two-step scheme. To simplify the explanation, 
assume that the first-component Markov chain $\{X^t\}$ is
regenerative with regeneration times $0\le T_0 < T_1 <\ldots$, so that
all values $X^{T_n}$ are equal to, say, $x^0=const$.

{\it Step 1}. Consider an embedded 2-component Markov chain
$(\widetilde{X}^n,\widetilde{Y}^n)=(X^{T_n},Y^{T_n})$
at the regenerative epochs $T_0<T_1<\ldots < T_n < \ldots$. 
Since $\widetilde{X}^{n}\equiv x^0$,
the sequence $\{\widetilde{Y}^n\}$ also forms a Markov chain.
For this chain one may show that, under appropriate assumptions and for an
appropriate test (Lyapunov) function $L$,
the drift ${\mathbf E} L(\widetilde{Y}^1 \ | \ \widetilde{Y^0}=y) - L(y)$
is bounded from above by the same constant for all values of $y$ and is also
uniformly negative if $y$ is outside of a certain
``bounded'' set, and hence this bounded set is positive recurrent.
For that, it is necessary to find or estimate from above
the average drift of $L(Y^n)$ over
a typical cycle, say $T_1-T_0$.
By introducing further smoothness condition(s), one can then ensure
that the Markov chain $\{\widetilde{Y}^n\}$ is also Harris ergodic.

{\it Step 2}. According to step 1, one can say that the
Markov chain $\{(X^n,Y^n)\}$ is regenerative, then verify its
aperiodicity and conclude that it is
Harris ergodic.

However, the first step of the proposed scheme may not be
implementable
if the transition function is not sufficiently smooth
(and, in particular, has discontinuities) as then it may be
difficult to find/estimate the value of the average drift of $L(Y)$ during a
typical cycle.

Another approach to proving stability is the use of fluid limits (see, e.g., \cite{dai1995_1}, \cite{dai1995_2} and Chapter 10 of \cite{meyn2007}). To apply this method, one needs to prove the convergence of an appropriately scaled version of the process, establish the dynamics of the limit, show that these dynamics are in some sense stable and then prove that this ``stability'' of the limit implies stability of the original Markov Chain. This is however often rather cumbersome and technically involved and moreover requires additional assumptions on the Markov chain.

In this paper we introduce a different approach
which is based on the following idea.

First, we find the (unique) stationary distribution $\pi_X$ for the
$X$-component. Second, we introduce an
auxiliary (time-homogeneous)
Markov chain $\{\widehat{Y}^t\}$ with transition probabilities
\begin{equation}\label{Yhat}
{\bf P} (\widehat{Y}^{t+1} \in \cdot \ | \widehat{Y}^t=y) =
\int_{\mathcal X} \pi_X (dx) {\mathbf P} (Y^1 \in \cdot
\ \big| \ X^1 = x, Y^0=y) \quad \mbox{a.s.}
\end{equation}
This Markov chain can also be viewed as an
outcome of the following recursive construction
of a two-dimensional Markov chain $(\widehat{X}^t,
\widehat{Y}^t)$ whose first components $\widehat{X}^t$ have
common distribution $\pi_X$.
First, for each $t$, we determine $\widehat{X}^{t+1}$
as a random variable which does not depend on all r.v.'s
$\{\widehat{X}^k\}_{k\le t},\{\widehat{Y}^k\}_{k\le t}$,
and has distribution $\pi_X$. Second, we determine
$\widehat{Y}^{t+1}$ as a random variable which is \\
(a) conditionally
independent of $\{\widehat{X}^k\}_{k\le t},\{\widehat{Y}^k\}_{k\le t-1}$
given $\widehat{X}^{t+1}$ and
$\widehat{Y}^t$, and \\
(b)
has distribution
\begin{eqnarray*}
 {\bf P} (\widehat{Y}^{t+1} \in \cdot \ |
\widehat{X}^{t+1}=x,\widehat{Y}^t=y) =
{\mathbf P} (Y^1 \in \cdot
\ \big| \ X^1 = x, Y^0=y) \quad \mbox{a.s.}
\end{eqnarray*}
Then, indeed, $\{\widehat{Y}^t\}$ forms a (time-homogeneous) Markov chain
with transition probabilities given by \eqref{Yhat}.

In this paper, we formulate and prove general stability criteria
in terms of averaging Lyapunov functions. In particular, we introduce 
natural conditions for the following implication to hold: if
the Markov chain $\{X^t\}$ is Harris ergodic then the Markov chain 
$\{\widehat{Y}^t\}$ is positive recurrent. The same conditions are then proved to imply that the
Markov chain $\{(X^t,Y^t)\}$ is positive recurrent too. 

We formulate and prove general stability results (Theorems 1 and 2) in Section 2.
In Theorem 1 we obtain a ``one-dimensional'' result and then in Theorem 2 its
multivariate analogue.
Then in Section~3 we study the stability of two systems with multiple access random
protocols in a changing environment. The proofs for their
stability are carried out by applying Theorem 2
and using the monotonicity arguments. In these particular scenarios we also show that
the stability conditions are necessary (if they do not hold, then the
system under
consideration is unstable).

\section{Stability of  two-component Markov chains using averaging
Lyapunov functions}

In this section we consider a general framework for the stability of a Markov Chain
containing
two components one of which is a Markov Chain itself.

In what follows, we write for short ${\mathbf P}_{x} (\ldots )$ instead of
${\mathbf P} (\ldots \ | \ X_0=x)$, ${\mathbf P}_{y} (\ldots )$ instead of
${\mathbf P} (\ldots \ | \ Y_0=y)$, and
${\mathbf P}_{x,y} (\ldots )$ instead of
${\mathbf P} (\ldots \ | \ X_0=x, Y_0=y)$. We use similar notation
${\mathbf E}_x$, ${\mathbf E}_y$, and ${\mathbf E}_{x,y}$
for conditional expectations.

\subsection{One-dimensional case}

Let $\{X^t\}$ and $\{Y^t\}$ be random sequences taking values in measurable spaces
$\left(\mathcal{X}, \mathcal{B_X}\right)$
and $\left(\mathcal{Y}, \mathcal{B_Y}\right)$, respectively, with countably
generated sigma-algebras ${\mathcal B_X}$
and ${\mathcal B_Y}$. Assume that $\{(X^t, Y^t)\}$ is
a Markov Chain on the state space
$\left(\mathcal{X} \times \mathcal{Y}, \mathcal{B_X} \times \mathcal{B_Y} \right)$ (here
$\mathcal{B_X} \times \mathcal{B_Y}$ is the minimal sigma-algebra
generated by sets
$B_1\times B_2$ where $B_1\in \mathcal{B_X}$ and
$B_2\in \mathcal{B_Y}$).
 Assume also
 the following:

$A0. \quad$ $\{X^t\}$ is a Markov Chain 
with ``autonomous'' dynamics:
for all $x\in \mathcal{X}$ and $y\in \mathcal{Y}$,
$$
{\mathbf P}_{x,y} (X^1 \in \cdot ) =
{\mathbf P}_x (X^1 \in  \cdot ).
$$
Further, the Markov chain $\{X^t\}$ satisfies the following conditions:

$A1. \quad$ It is aperiodic;

$A2. \quad$ There exists a set~$V \in \mathcal{B_X}$ such that
\begin{equation} \label{eq:def_recurrence_1}
\tau \equiv \tau (V) = \min\{t \ge 1: X^t \in V\} < \infty \quad
{\mathbf P}_x-\text{a.s.}
\end{equation}
for any initial value $X_0=x \in \mathcal{X}$
and, moreover,
\begin{equation} \label{eq:def_recurrence_2}
s_0=\sup_{x \in V} \mathbf E_x \tau  < \infty.
\end{equation}
(We say that the set~$V$ is {\it positive recurrent}
 for~$\{X^t\}$ if conditions~(\ref{eq:def_recurrence_1}) and~(\ref{eq:def_recurrence_2}) hold.)

$A3. \quad$ The set $V$ admits a minorant measure, or is petite (in the terminology
of \cite{MeynTweedie}), i.e. there exist a number~$0 < p \le 1$, a measure~$\mu$ and an
integer~$m \ge 1$ such that
$$
{\mathbf P}_x (X^m \in \cdotp ) \ge p \mu(\cdotp)
$$
for any~$x \in V$.

Conditions A1--A3 imply that the $X$-chain is {\it Harris ergodic}, so there exists a stationary
distribution~$\pi = \pi_X$ such that
\begin{equation}\label{uniformx}
\sup_{B \in \mathcal{B}} \left|{\mathbf P}_x \left(X^t \in B \right) - \pi(B)\right| \to 0
\end{equation}
as~$t \to \infty$, for any~$x \in \mathcal{X}$. 
Moreover, Conditions A1-A3 imply that
\begin{equation}\label{uniformss}
\mbox{convergence in} \ \  \eqref{uniformx} \ \
\mbox{is uniform in} \ \ x\in V,
\end{equation}
see e.g.
\cite{Borovkov}, \cite{Thorisson}
or \cite{Lindvall}. We can formulate also a coupling version of conditions
\eqref{uniformx}--\eqref{uniformss}:\\
for any $x\in V$, there is a coupling of $\{X^{t}\}$ and of a
stationary Markov chain $\{ \widehat{X}^{t}\}$ having distribution
$\pi$ such that
\begin{equation}\label{c_uni1}
  \nu = \min \{t \ : \ X^{k}=\widehat{X}^{k}, \  \mbox{for all}
\ \ k\ge t\} < \infty \ \ {\mathbf P}_x-\mbox{a.s. }
\end{equation}
for $X_{0}=x$ and
\begin{equation}\label{c_uni2}
\delta_{t} :=\sup_{x\in V} {\mathbf P}_x (\nu >t) \to 0, \quad t\to\infty,
\end{equation}
see Appendix A for a proof of~\eqref{c_uni1} and~\eqref{c_uni2}.

Introduce a function
\begin{eqnarray} \label{eq:def_lx}
L_1(x) =
\begin{cases}
\mathbf E_x \tau, \quad \text{if} \quad x \notin V, \\
0, \quad \text{if} \quad x\in V.
\end{cases}
\end{eqnarray}
Then
\begin{equation} \label{eq:lyapunov_x1}
{\mathbf E}_x \left(L_1(X^1) - L_1(x)\right) = -1
\end{equation}
for all $x \notin V$ and
\begin{equation} \label{eq:lyapunov_x2}
{\mathbf E}_x \left(L_1(X^1) - L_1(x)\right) \le \sup_{x \in V} \mathbf E_x\tau -1
< \infty
\end{equation}
for all $x \in V$. These inequalities follow from observing that
$\tau_{X_1} = \tau_{x} - 1$ ${\mathbf P}_x-\mbox{a.s.}$  
if $X_1 \notin V$.
Inequalities~(\ref{eq:lyapunov_x1}) and~(\ref{eq:lyapunov_x2}) mean that the function~$L_1$ is
an appropriate Lyapunov function 
for the Markov Chain~$\{X^t\}$ in the sense that it satisfies the standard conditions
for the Foster criterion to hold.

It is known (see, e.g., \cite{Borovkov}, \cite{Thorisson} or \cite{Lindvall})  that Conditions A1--A3 imply that
\begin{equation}\label{olittle}
\lim_{t\to\infty}\frac{1}{t} \sup_{x\in V} {\mathbf E}_x L_1(X^t) =0.
\end{equation}

$B. \quad$ For the sequence~$\{Y^t\}$ we assume that
there exists a non-negative measurable function~$L_2$
such that:

$B1. \quad$ The expectations of the absolute values of the
increments of the sequence $\{ L_2(Y^t)\}$ are bounded
from above by a constant $U$:
$$
\sup_{x \in \mathcal{X}, y \in \mathcal{Y}} {\mathbf E}_{x,y} \left|L_2\left(Y^{1}\right) - L_2\left(Y^0\right)\right| \le U
<\infty. 
$$

$B2. \quad$ 
There exist a non-negative and non-increasing function $h(N), N\ge 0$ such
that
$h(N) \downarrow 0$
as $N\to\infty$, and a measurable function $f : \mathcal{X} \to
(-\infty ,\infty )$ such that
$$
\int_{\mathcal{X}} f(x) \pi (dx) := -\varepsilon < 0
$$
and
\begin{equation}\label{mainineq}
{\mathbf E}_{x,y}\left(L_2\left(Y^1\right) - L_2(y)\right) \le f(x) +
h(L_2(y))
\end{equation}
for all $x\in \mathcal{X}$ and $y\in\mathcal{Y}$.

It follows from Condition B1 
that, 
without loss of generality,
the function $f$
may be assumed to be bounded:
\begin{equation}\label{boundedf}
\sup_{x\in\mathcal{X}}|f(x)|=K<\infty.
\end{equation}

\begin{remark}
Note that conditions B are sufficient for stability of the chain $\widehat{Y}$ introduced in (\ref{Yhat}), due to the standard Foster criterion. We are going to show that they (together 
with conditions A) are also sufficient for stability of the original chain $\{Z^t\}$.
\end{remark}

\begin{lemma}\label{base1}
Conditions A and B imply that, for some positive integer
$t_0$ and for any integer~$t\ge t_0$ there exists a positive number~$N_0 =
N_0(t_0)$
such that 
if $L_2(y)\ge N_0$, then
\begin{equation}\label{power}
{\mathbf E}_{x,y} \left( L_2(Y^{t}) -
L_2(y)\right) \le - t \Delta
\end{equation}
for all $x\in V$ with~$\Delta = \varepsilon/10$.
\end{lemma}

{\sc Proof.} Due to conditions \eqref{uniformss} and \eqref{boundedf},
for any $c_1\in (0,1)$ one can choose a number $n_0$ such that
\begin{equation}\label{firstbound}
\sup_{t\ge n_0} \sup_{x\in V} \left|
{\mathbf E}_x \left(f\left(X^t\right)\right) -
\int_{\mathcal{X}} f(z) \pi (dz) \right| \le c_1.
\end{equation}
Let $t=n_0+m$. Then, for $x\in V$, $y\in\mathcal{Y}$ and for $U$ from Condition B1,
\begin{eqnarray*}
{\mathbf E}_{x,y}
( L_2(Y^{t}) - L_2(y))
&=&
\sum_{i=0}^{t-1} {\mathbf E}_{x,y} \left( L_2\left(Y^{i+1}\right)
- L_2\left( Y^i\right)\right)\\
&\le & Un_0 +\sum_{i=n_0}^{t-1}
{\mathbf E}_{x,y}
\left(
{\mathbf E}
\left(
L_2\left(Y^{i+1}\right) \ | \ X^i,Y^i \right)
- L_2\left(Y^i\right)
\right)\\
&\le &
Un_0 + \sum_{i=n_0}^{t-1} {\mathbf E}_{x,y}
\left(
f\left(X^i\right)+h\left(Y^i\right) \right)\\
&\le &
Un_0 - (\varepsilon - c_1) m
+\sum_{i=n_0}^{t-1} {\mathbf E}_{x,y} h\left(Y^i\right)\\
&\le &
Un_0 - m \left(\varepsilon - c_1
- h(0) {\mathbf P}_{x,y}\left(\min_{i\le t} L_2(Y^i) < \widehat{N}\right)
- h(\widehat{N})\right)
\end{eqnarray*}
where $\widehat{N}$ is any positive number.
Take $c_1=\varepsilon /5$ and then $m_0 = \max (n_0, 5Un_0/\varepsilon )$
and $t_0=n_0+m_0$.
Then, for any $m\ge m_0$, let $\widehat{N}$ be such that
$h(\widehat{N})\le \varepsilon /5$ and $h(0) (n_0+m)^2 U/\widehat{N} \le \varepsilon /5$.
Further, let $N_0=2\widehat{N}$. If $L_2(y)\ge N_0$, then
\begin{eqnarray*}
{\mathbf P}_{x,y}\left(\min_{1\le i\le t} L_2(Y^i) < \widehat{N}
\right)
&\le & t \sup_{z,u}{\mathbf P}_{z,u} \left(
L_2(Y^1)-L_2(u) \le -\widehat{N}/t \right)\\
&\le &
h(0)t^2 U/\widehat{N}
\le
\varepsilon /5,
\end{eqnarray*}
where we applied the Markov inequality.
Hence, for $x\in V$ and $y$ such that $L_2(y)\ge N_0$,
$$
{\mathbf E}_{x,y}
( L_2(Y^{t}) - L_2(y)) \le -m \frac{\varepsilon}{5} \le - t \frac{\varepsilon}{10}
$$
and therefore inequality \eqref{power} holds with $\Delta = \varepsilon /10$.
\qed

\begin{theorem} \label{th:stability_general}
Under the conditions~$A$ and~$B$ there exists $N_0$ such that the set
$D:=V \times \{y: L_2(y) \le N_0\}$ is {\it positive recurrent}
for the Markov Chain~$\{(X^t, Y^t)\}$. 
\end{theorem}

\begin{cor}
Assume Conditions A and Condition B1 to hold. Assume further
that,
uniformly in $x\in \mathcal{X}$, the conditional distribution
$
{\mathbf P}_{x,y} (L_2\left(Y^1\right) - L_2(y) \in \cdot )
$
converges weakly and in  $\mathcal{L}_1$ to a limiting one,
say $H_x(\cdot )$, as $L_2(y)\to\infty$,
where
$$
\int_{\mathcal{X}} \int_{\mathcal{Y}} \pi (dx) y H_x(dy)
$$
is negative.
Then the Markov Chain~$\{(X^t, Y^t)\}$ is positive recurrent.
\end{cor}

{\sc Proof of Theorem~\ref{th:stability_general}.}
We use the following extension of the Foster
criterion (see, e.g., \cite[Theorem 1]{Foss} or \cite{Bramson}
and references therein).
Let $L : \mathcal{X}\times\mathcal{Y}\to [0,\infty )$ be a measurable
function. Then a set $D=\{ (x,y) \ : \ L(x,y)\le N\}$ is positive
recurrent if there exists a positive integer-valued measurable
function $T$ on  $\mathcal{X}\times\mathcal{Y}$ such that
$$
\sup_{x,y} \frac{T(x,y)}{\max (1,L(x,y))} <\infty ,
$$
$$
\sup_{(x,y)\in D} {\mathbf E}_{x,y} L\left(X^{T(x,y)},Y^{T(x,y)}\right)
< \infty
$$
and
\begin{equation}\label{Txy}
{\mathbf E}_{x,y} L\left(X^{T(x,y)},Y^{T(x,y)}\right) -
L(x,y) \le -cT(x,y)
\end{equation}
for some $c>0$ and for all $(x,y)\notin D$.

Let $H > U$ with $U$ from condition $B1$ and
let $t_0$ be such that
\begin{equation}\label{olittle2}
\sup_{t\ge t_0} \frac{1}{t} \sup_{x\in V} {\mathbf E}_x L_1(X^t) \le \frac{\Delta}{2H},
\end{equation}
which is possible due to~\eqref{olittle}. Here again $\Delta = \varepsilon /10$.
Let $n_0$,$m$ and $\widehat{N}$ be chosen according to Lemma 1 with $t_0$ from
\eqref{olittle2}. Let again $N_0=2\widehat{N}$.

We take the set $D=V \times \{y: L_2(y) \le N_0\}$ 
and define a function
$$
L(x,y) = H L_1(x) + L_2(y).
$$
Further, take $T(x,y)=1$ if either $(x,y)\in D$ or $x\notin V$ and
$T(x,y)=t$ if $x\in V$ and $L_2(y)> N_0$, where $t=n_0+m$ is from the
proof of Lemma 1.

Now, if  $(x,y)\in D$, then
$$
{\mathbf E}_{x,y} L(X^1,Y^1) - L(x,y) \le H \sup_{x\in V} {\mathbf E}_x \tau + U.
$$
If $x\notin V$, then
$$
{\mathbf E}_{x,y} L(X^1,Y^1) - L(x,y) \le - H+ U <0.
$$
Finally, if  $x\in V$ and $L_2(y)> N_0$, then, by \eqref{power},
$$
{\mathbf E}_{x,y} L\left(X^{T(x,y)},Y^{T(x,y)}\right) -
L(x,y) \le H t \frac{\Delta}{2H} - t \Delta = -t \frac{\Delta}{2}.
$$
\qed

\subsection{Multivariate case}

Now we formulate and prove a multivariate analogue of Theorem 1.
We modify the model as follows. We continue to assume that $\{(X^t, Y^t)\}$ is
a Markov Chain on the state space
$\left(\mathcal{X} \times \mathcal{Y}, \mathcal{B_X} \times \mathcal{B_Y} \right)$
and that $\{ X^t\}$ is a Markov chain satisfying conditions A0--A3. But now we
assume that the state space $\mathcal{Y}$ is a product of $M$ spaces
$\mathcal{Y}= \widetilde{\mathcal{Y}_1}\times \ldots \times \widetilde{\mathcal{Y}_M}$
(with the product sigma-algebra), so the $Y$-component of the Markov chain
has $M$ coordinates, $Y^t = (Y^t_1,\ldots ,Y^t_M)$. Clearly, the results of the previous section hold in this case too. However, they turn out not to be applicable in the examples we are going to consider in the second part of the paper, and in this section we develop conditions only involving each individual coordinate of the $Y$-chain. We believe that these conditions may be useful in many applied models and, in particular, where the $Y$-chain has a number of coordinates which 
may be dependent but have simple individual dynamics.

We assume that there is a non-negative
function $L_{2,i}$  defined on $\widetilde{\mathcal{Y}_i}$ and that conditions similar to
$B1-B2$ hold for each coordinate. 

$\widetilde{B1}. \quad$ For any $i=1,\ldots ,M$, the expectations
of the absolute values of the
increments of the sequence
$\{L_{2,i}(Y^t_i)\}$ are bounded from above by a constant $U$:
$$
\sup_{x \in \mathcal{X}, y \in \mathcal{Y}}{\mathbf E}_{x,y} \left|L_{2,i}\left(Y^1_i\right) -
L_{2,i}\left( Y^0_i\right)\right| \le U < \infty.
$$

$\widetilde{B2}.\quad$
For each~$i$, there exists a function $h_i(N), N\ge 0$ such
that $h_i(N) \downarrow 0$
as $N\to\infty$ and a measurable function $f_i : \mathcal{X} \to
(-\infty ,\infty )$ such that $\sup_x |f_i(x)| :=K_i < \infty$,
\begin{equation}\label{fdrift}
\int_{\mathcal{X}} f_i(x) \pi (dx) := -\varepsilon_i < 0
\end{equation}
and
\begin{equation}\label{mainineq_tilde}
{\mathbf E}_{x,y}\left(L_{2,i}\left(Y_i^1\right) - L_{2,i}(y_i)\right) \le f_i(x) +
h_i(L_{2,i}(y_i))
\end{equation}
for all $x\in \mathcal{X}$ and $y\in\mathcal{Y}$, where the RHS depends only on the $x$ and on the $i$th coordinate of the $y$.

In the multivariate case, the
statement of
Lemma~\ref{base1}
holds for each coordinate $i$, and we restate it here for convenience.

\begin{lemma} \label{base2}
Conditions $A$ 
and $\widetilde{B1}$--$\widetilde{B2}$ imply that,  
for some
integer~$t'$ and for any integer $t \ge t'$,
there exists a positive number~$N_0$ such that, for any $i$,
if $L_{2,i}(y)\ge N_0$ then for all $x\in V$
\begin{equation}\label{power2}
\mathbf{E}_{x,y} \left( L_{2,i}(Y^{t}) -
L_{2,i}(y)\right) \le - t \Delta,
\end{equation}
where $\Delta = \min \varepsilon_i/10$.
\end{lemma}

If the sequence~$\{(X_1^t,..,X_M^t, Y_i^t)\}$ is a Markov chain for each~$i$, then clearly, Theorem 1 also holds for each coordinate $i$, meaning that the Markov chain~$\{(X_1^t,..,X_M^t, Y_i^t)\}$ for any coordinate $i$ is positive recurrent. However, this does not imply the positive recurrence of the entire chain \linebreak $\{(X_1^t,..,X_M^t, Y_1^t,..,Y_M^t)\}$. 
In order to obtain conditions for stability in the multivariate case (also without assuming that each of the above sequences is itself a Markov chain), we
require an extra assumption to hold:

$\widetilde{B3}. \quad$ For each~$i$ it holds that
$$
\sup_{t \ge N} \sup_{y \in \mathcal{Y}, x \in V} \frac{\mathbf{E}_{x,y} \left(L_{2,i}(Y^t) - L_{2,i}(y)|Y^0 = y \right)}{t} \to 0
$$
as~$N \to \infty$.

Assumption~$\widetilde{B3}$ means that the drift of the
function~$L_{2,i}$ in~$t$ steps grows slower than any linear
function of~$t$. 
Below we provide a simple condition $\widehat{B}$ which is sufficient
for $\widetilde{B3}$ to hold. Note that condition $\widehat{B}$ is
 a
 stronger version
of condition $\widetilde{B}_2$ and that
it can be easily verified in many applications.

\begin{lemma}\label{suffi}
Assume that the following condition holds:\\
$\widehat{B}. \quad$ For each $i$, there exist
a non-negative and non-increasing function $h_i(N)$,$N\ge 0$ such that
$h_i(N)\downarrow 0$ as $N\to\infty$,
and a family of mutually independent
random variables $\{ \varphi_{x,i}^t \},  x\in \mathcal{X}, t=0,1,\ldots$
such that,
\begin{itemize}
  \item[(i)] for each~$t$ and~$i$ these random variables are uniformly integrable;
  \item[(ii)]
for each $i$ and $x$, the random variables
$\{\varphi_{x,i}^t, \ t=0,1,\ldots \}$ are identically distributed with common distribution
function $F_{x,i}$, which is such that $F_{x,i}(y)$ is measurable as
a function of $x$ for any fixed $y$;
  \item[(iii)] the inequality
\begin{equation}\label{multiX}
L_{2,i}\left(Y_i^{t+1}\right) - L_{2,i} \left(Y_i^t\right)
\le \varphi_{X^t,i}^t + h_i(L_{2,i}(Y_i^t)) \quad \mbox{a.s.}
\end{equation}
holds for all $x\in \mathcal{X}$,
$y\in \mathcal{Y}$ and $t=0,1,\ldots$;
  \item[(iv)]
functions
$f_i(x) = {\mathbf E} \varphi_{x,i}^1$ satisfy
 condition
\eqref{fdrift}.
\end{itemize}
Then assumption~$\widetilde{B3}$ holds too.
\end{lemma}

We will now formulate the main theorem and then prove it. 
A proof of Lemma \ref{suffi} is presented after the proof
of the main theorem.

\begin{theorem} \label{th:stability_multi}
Under assumptions~$A$ and
assumptions~$\widetilde{B1}-\widetilde{B3}$ (or
assumptions~$\widetilde{B1}$ and $\widehat{B}$),
there exists $N_1\ge N_0$ such that the set
$D:=V \times \{y: \sum\limits_{1=i}^M L_2(y_i) \le N_1\}$ is {\it positive recurrent}
for the Markov Chain~$\{(X^t, Y^t)\}$. 
\end{theorem}

{\sc Proof of Theorem 2.} We may apply Lemma 2 to each of~$Y_i$ and may therefore assume that
inequality \eqref{power2} holds for each coordinate, with the same~$t'$ and
the same~$\Delta$.

Due to condition~(\ref{olittle}), we may also assume that~$t'$ is such that
\begin{equation}\label{olittle3}
\sup_{t\ge t'} \frac{1}{t} \sup_{x\in V} {\mathbf E}_x L_1(X^t) \le \frac{\Delta}{2H},
\end{equation}
where~$H$ is any positive number larger than~$MU$ with~$U$ from assumption $\widetilde{B1}$.

Choose also~$n_0$ such that
\begin{equation} \label{eq:sublinear}
\mathbf{E}_{x,y} \left(L_{2,i}(Y^t) - L_{2,i}(y)|Y^0 = y \right) \le \frac{\Delta}{2M} t
\end{equation}
for all~$t \ge n_0$ and for all~$y$. This is possible due to condition~$\widetilde{B3}$.

We again use Theorem~1 of~\cite{Foss}.
Take~$N_1 = MN_0$ so that~$D = V \times \{y: \sum\limits_{i=1}^M L_2(y_i) \le M N_0\}$ and introduce test function
$$
L(x,y) = HL_1(x) + \sum_i L_{2,i}(y_i).
$$
Let $T(x,y)=1$ if either $(x,y)\in D$ or $x\notin V$ and
$T(x,y)= t_1 := \max\{t', n_0\}$, otherwise.

Now, if  $(x,y)\in D$, then
$$
{\mathbf E}_{x,y} L(X^1,Y^1) - L(x,y) \le H \sup_{x\in V} {\mathbf E}_x \tau + MU < \infty.
$$
If $x\notin V$, then
$$
{\mathbf E}_{x,y} L(X^1,Y^1) - L(x,y) \le - H+ MU <0.
$$
And if  $x\in V$ and $\sum\limits_{i=1}^M L_{2,i}(y_i)> MN_0$, then $L_{2,i}(y_i)>N_0$ for at least one index~$i$. Denote by~$k$ the number of such indices. Then,
by \eqref{power2}, \eqref{olittle3} and \eqref{eq:sublinear},
\begin{multline*}
{\mathbf E}_{x,y} L\left(X^{T(x,y)},Y^{T(x,y)}\right) -
L(x,y) \le H t_1 \frac{\Delta}{2H} - k t_1 \Delta + (M-k) \frac{\Delta}{2M} t_1 \le - t_1 \frac{\Delta}{2M}
\end{multline*}
as~$k \ge 1$. 
Clearly,
$
\sup_{x,y} \frac{T(x,y)}{\max (1,L(x,y))} <\infty,
$
so the theorem is proved.
\qed

{\sc Proof of lemma \ref{suffi}}.
Let $C_1 = \max_{1\le i \le M} h_i(0)$ and $C_2>0$ be such that
$h_i(C_2) \le \varepsilon_i/2$ for all $i$.
Then inequality \eqref{multiX} implies that
\begin{equation}\label{newineq}
L_{2,i}(Y_i^{t+1})-L_{2,i}(Y_i^t)\le
\psi_{X^t,i}^t+C_1{\bf I}(L_{2,i}(Y_i^t)\le C_2) \quad
\mbox{a.s.,}
\end{equation}
where $\psi_{X^t,i}^t=\varphi_{X^t,i}^t+\varepsilon_i/2.$

Let $\nu_t\le t$ be the last time $k$ before $t$ when
$L_{2,i}(Y_i^k)\le C_2$ (we let $\nu_t=0$ if such $k$ does not
exist). Then
\begin{eqnarray*}
L_{2,i}(Y_i^{t+1})-L_{2,i}(Y_i^0) & \le &
C_1+C_2 + \sum_{k=\nu_t}^t \psi_{X^{k},i}^k \\
&\le &
C_1+C_2 + \max_{0\le j \le t} \sum_{k=j}^t
\psi_{X^{k},i}^k.
\end{eqnarray*}
For each $x\in V$, let $\widehat{X}^{t}$ be a stationary sequence
satisfying conditions~\eqref{c_uni1}-\eqref{c_uni2}. Then
\begin{eqnarray*}
  \frac{1}{t}\max_{0\le j \le t} \sum_{k=j}^t
\psi_{X^{k},i}^k &\le &
\frac{1}{t} \max_{0\le j \le t} \sum_{k=j}^t
\psi_{\widehat{X}^{k},i}^k +
\frac{1}{t} \sum_0^{t} \left|
 \psi_{{X}^{k},i}^k -
 \psi_{\widehat{X}^{k},i}^k
 \right|.
   \end{eqnarray*}
Here
the first summand on the RHS tends to 0 both a.s. and in mean,
since the sequence $\{ \psi_{\widehat{X}^{k},i}^k \}$
is stationary ergodic with negative mean
${\mathbf E} \psi_{\widehat{X}^1,i}^1 \le -\varepsilon_i/2$.
The second summand also tends to 0
in mean uniformly in $x\in V$, due to uniform integrability.
Indeed, for any $\Delta >0$ let $R$ be such that
${\mathbf E} \left|\psi_{x,i}^{1} {\bf I}(|\psi_{x,i}^{1}| > R)\right|
\le \Delta$, for all $x\in X$. Then for any $x\in V$ and any $k$
\begin{eqnarray*}
{\mathbf E} \left|
  \psi_{X^{k},i}^{k} \right|{\bf I} \left(
  \psi_{X^{k},i}^{k}\ne \psi_{\widehat{X}^{k},i}^{k}\right) \
 &\le &
 \Delta + R\delta_{k}
\end{eqnarray*}
and a similar inequality holds for
$\psi_{\widehat{X}^{k},i}^{k}$. Therefore,
\begin{eqnarray*}
{\mathbf E}
\frac{1}{t} \sum_0^{t} \left|
 \psi_{{X}^{k},i}^k -
 \psi_{\widehat{X}^{k},i}^k
 \right| &\le &
 \frac{1}{t}\sum_0^t{\mathbf E}
 \left(
 \left| \psi_{X^k,i}^k\right|
  +
 \left| \psi_{\widehat{X}^k,i}^k\right| \right){\bf I}
 \left(
  \psi_{\widehat{X}^{k},i}^{k}\ne \psi_{\widehat{X}^{k},i}^{k}\right)\\
 &\le &
   2\Delta + 2R \frac{1}{t} \sum_0^t\delta_{k}.
   \end{eqnarray*}
   As $\Delta$ is arbitrary small and
   $\delta_{t}\to 0$, the second summand on the RHS of the latter inequality may be
   made arbitrarily small, and the result follows.
\qed

\section{Stability of two systems with two random multiple-access protocols and a
finite number of stations}

In order to improve performance and usability of laptops, tablet computers and other portable equipment, such devices are nowadays often endowed with several wireless interfaces. For example, if a device is equipped with the WiFi interface (protocol standard IEEE~802.11) and the WiMax interface (protocol standard IEEE~802.16), the user has a wide range of possibilities for the wireless exchange of data. However, it is known that if the same device uses both these standards, their data streams interfere with each other and the speed of transmission and throughputs therefore deteriorate. There are currently attempts to minimise the effect of such interference (see, e.g. \cite{Andreev}, \cite{Zhu}).

In this paper we introduce a simplified model that mimics the essential features of both standards and of the interference between them described in the previous paragraph. The IEEE~802.11 standard uses a random-access algorithm for message transmissions, therefore multiple message transmissions may be attempted simultaneously and collide with each other. Such conflicts are then resolved according to a specific algorithm. In the IEEE~802.16 standard message transmission follow a schedule and no conflicts occur. In our model, we assume that in the random-access algorithm every station transmits with a certain probability, independently of everything else (cf. ALOHA algorithm, see~\cite{Abramson} and~\cite{Roberts}). Each message that experienced a collision then simply returns to its origin and is treated as any other message. In other words, no conflict-resolution mechanism is implemented. Our model also assumes that the scheduling is completely symmetric: every node receives time intervals of the same duration to transmit their messages, successively.

We consider two cases. In the first case the
system has the so-called MAC-coordinator (see~\cite{Andreev} and~\cite{Zhu}) which gives a priority to the IEEE~802.16 protocol. This case is considered in subsection~\ref{subsec:blocking}. The second case concerns a system without a MAC-coordinator and is considered in subsection~\ref{subsec:no_blocking}.

\subsection{Network with a MAC-coordinator} \label{subsec:blocking}

Assume there are~$M$ identical stations numbered~$1,...,M$. There are~$2$ types of messages called
``red" and ``green", and each station has two infinite buffers where these messages may be stored, one for each type.

We make a few assumptions:

{\bf Assumption 1.} Time is slotted, stations may only start transmissions at the beginning of a time slot,
and each transmission time is equal to the length of a slot. Hence, we may assume that the events (such
as arrival of a new message, the beginning of a transmission and the end of a transmission) may only happen
at time instants $1, 2, ....$. We also assume that the transmission channel is such that, at a given time
slot, at most one red and at most one green message may be transmitted. Note also that any single station
cannot transmit two messages (red and green) simultaneously during the same time slot. Summarising, in
any time slot, there may be no transmissions at all. Otherwise,
 there may be a transmission of either only one red message or only
one green message or
one red and one green message, but in the latter case the transmissions have to be made by different stations.

{\bf Assumption 2.}  Transmissions of red messages are scheduled and do not collide.
More precisely, for $t=1,2,\ldots$,
time slot~$t$ is scheduled for a transmission from node~$i(t) = ((t-1) \mod M) + 1$:
if the queue of red messages at that node is non-empty, then there is a (successful)
transmission of the first one of them; otherwise there is no transmission of red
messages at that time slot.

{\bf Assumption 3.} Transmissions of green messages follow the well-known ALOHA protocol (see \cite{Roberts}):
at a given time slot,
every node that is not transmitting a red message and whose queue
of green messages is not empty, transmits a green message (say, first in the queue)
with probability $p$, independently of everything else.
Then one of three events may
occur:
\begin{itemize}

\item{Only one node attempts to transmit a green message.
Then the transmission is successful;}

\item{No transmission attempted;}

\item{Two or more nodes attempt transmissions of green messages.
Then all these transmissions fail due to collision, and the messages
stay in their queues.}

\end{itemize}


{\bf Assumption 4.} Red messages arrive in the system as an i.i.d.
sequence $\{ \xi^t\}_{t\ge 0}$ with a finite
rate~$\lambda_R = \mathbf{E} \xi^t$ 
(here $\xi^t$ is the total
number of new red messages within time slot $(t-1,t]$).
Similarly,  green
messages arrive independently as a renewal sequence $\{
\eta^t\}_{t\ge 0}$ with a finite rate~$\lambda_G = \mathbf{E} \eta^t$.
Every
arriving message is assigned to a node at random with equal
probabilities~$1/M$.


Let~$R_i^t$ and~$G_i^t$ be the numbers of red and green messages respectively in the queue
of node $i$ at
the beginning of a time slot $t$. The sequence $\{(R_1^t,...,R_M^t)\}$ forms
 a Markov Chain, and so does the sequence $\{(R_1^t,...,R_M^t,G_1^t,...,G_M^t)\}$. Note also
 that the latter Markov Chain describes the state of the system completely. We will say that the
 system is {\it stable} if its underlying Markov Chain is positive recurrent.

For such an algorithm we prove the following

\begin{theorem} \label{th:aloha}
Assume $\lambda_R < 1$. Then the system is stable if
\begin{equation} \label{eq:cond_aloha}
\begin{cases}
\lambda_G < (1-\lambda_R) p, \quad \text{if} \quad M=1, \cr
\lambda_G < \lambda_R (M-1) p (1-p)^{M-2} + (1-\lambda_R)
Mp(1-p)^{M-1}, \quad \text{if} \quad M>1
\end{cases}
\end{equation}
and unstable if the opposite strict inequality holds.
\end{theorem}

\begin{remark}
The conditions of the theorem seem to be intuitively clear. Consider a single time slot. The average number of new green messages in a time slot is equal to~$\lambda_G$. Let us now consider the average number of green messages leaving the system. For that, assume that all nodes have a green message (this may be thought as the ``worst-case" scenario). Assume also that the Markov chain representing the states of the queues of red messages in all nodes is in its stationary regime. In that case, when we look at a single time slot, the queue of red messages of the node that is scheduled to transmit a red message in the time slot will be non-empty with probability~$\lambda_R$ (this is quite clear intuitively but see Appendix B for a proof) and empty with probability~$1-\lambda_R$. In the former case, there are~$M-1$ nodes that can successfully transmit a green message and probability of success for each of them is~$p(1-p)^{M-2}$ and therefore the expected number of green messages leaving the system is~$(M-1)p(1-p)^{M-2}$. The latter case may be considered in a similar way.

The simple heuristic argument presented above depends heavily on the symmetry of the system. The symmetry is also essential for our strict analysis.
\end{remark}

{\sc Proof of Theorem~\ref{th:aloha}.}
We give a proof only for~$M > 1$. The case~$M=1$
may be considered following the same lines and applying straightforward 
changes. We start by proving stability, a proof of instability follows. 

{\sc Proof of stability}.
We apply Theorem~\ref{th:stability_multi}. Denote
by~$\xi_i^t$ and $\eta_i^t$ the numbers of new red and,
respectively, green packets arriving at time slot~$t$ to
station~$i$. Since the total number of red arrivals in time slot $t$ is
$\xi^t$ and each of them chooses one of the stations at random, each
$\xi_i^t$ has a conditional binomial distribution with parameters $\xi^t$ and
$1/M$ and $\sum_1^M \xi_i^t = \xi^t$. Similarly, for $t=0,1,\ldots$,
a random variable
$\eta_i^t$ has 
a binomial distribution with parameters $\eta^t$ and $1/M$. Denote also by
\begin{eqnarray*}
\alpha_i^t = \begin{cases}
1, \quad \text {with probability} \quad p, \\
0, \quad \text {with probability} \quad 1-p
\end{cases}
\end{eqnarray*}
the sequence of i.i.d. random variables representing the decisions taken by the nodes on
whether or not to attempt a transmission of a green message.
Then the Markov Chain~$\{(R_1^t,...,R_M^t,G_1^t,...,G_M^t)\}$ has the following transitions:
\begin{eqnarray*}
R_i^{t+1} = \begin{cases}
R_i^t + \xi_i^t, \quad \text{if} \quad i \neq i(t), \\
R_i^t + \xi_i^t - \I\{R_i^t > 0\}, \quad \text{if} \quad i = i(t),
\end{cases}
\end{eqnarray*}
where, as before,~$i(t) = ((t-1) \mod M) + 1$. Further,
let $\gamma_j^t = \alpha_j^t \I \{G_j^t>0\}$. Then
\begin{eqnarray*}
G_i^{t+1} = \begin{cases} G_i^t + \eta_i^t - \gamma_i^t
\prod\limits_{j \neq i, j \neq i(t)} \left(1-\gamma_j^t\right)
\left(1-\gamma_{i(t)}^t  \I\{R_{i(n)}^t = 0\}\right), \quad
\text{if} \quad i \neq i(t), \\
G_i^t + \eta_i^t - \I\{R_i^t = 0\} \gamma_i^t
\prod\limits_{j \neq i}\left(1-\gamma_j^t\right), \quad \text{if} \quad i=i(t).
\end{cases}
\end{eqnarray*}

Introduce now a new Markov Chain~$\{(R_1^t,...,R_M^t,\widetilde{G}_1^t,...,\widetilde{G}_M^t)\}$
where the first~$M$ components (representing the states of the red queues)
are the same as before,
and the remaining components (representing the states of the green queues) have the following transitions:
\begin{eqnarray*}
\widetilde{G}_i^{t+1} = \begin{cases} \widetilde{G}_i^t + \eta_i^t -
\widetilde{\gamma}_i^t
\prod\limits_{j \neq i, j \neq i(t)} \left(1-\alpha_j^t\right)
\left(1-\alpha_{i(t)}^t \I\{R_{i(t)}^t = 0\}\right), \quad \text{if} \quad i \neq i(t), \\
\widetilde{G}_i^t + \eta_i^t - \I\{R_i^t = 0\} \widetilde{\gamma}_i^t
\prod\limits_{j \neq i}\left(1-\alpha_j^t \right), \quad \text{if} \quad i=i(t).
\end{cases}
\end{eqnarray*}
Here $\widetilde{\gamma}_j^t = \alpha_j^t \I \{\widetilde{G}_j^t >0\}$, for all $j$ and $t$.

In words, the Markov Chain~$\{(R_1^t,...,R_M^t,\widetilde{G}_1^t,...,\widetilde{G}_M^t)\}$
represents the state of the system where each station with an empty green queue (if not blocked by
a transmission of a red message) may send (and does so with probability~$p$) a ``dummy" packet which
interferes with (dummy or legitimate)
packets of other stations.

Since~$\gamma_i^t \le \alpha_i^t$ a.s. for all~$i$ and~$t$, it follows from the two systems of equations displayed above that the new Markov Chain
dominates the initial one:
if~$\{(G_1^1,...,G_M^1)\} = \{(\widetilde{G}_1^1,...,\widetilde{G}_M^1)\}$, then
$$
\{(R_1^t,...,R_M^t,G_1^t,...,G_M^t)\} \le
\{(R_1^t,...,R_M^t,\widetilde{G}_1^t,...,\widetilde{G}_M^t)\}
$$
a.s., for any~$t$. Hence, to prove stability of the initial Markov Chain, it
is sufficient to prove stability of the new one. As it follows 
from Theorem~\ref{th:stability_multi}, Lemma~\ref{suffi} and the symmetry of the system, 
it is sufficient to show that condition~$\widetilde{B1}$ and conditions of Lemma~\ref{suffi} 
hold for the sequence~$\{(R_1^t,...,R_M^t,\widetilde{G}_1^t)\}$. Note that it forms a Markov chain in this case.
For simplicity, we will also omit the coordinate index~$1$ in all the functions appearing in conditions. Take function~$L_2(y) = y$ and consider the state of the Markov Chain~$\{(R_1^t,...,R_M^t,\widetilde{G}_1^t)\}$
after every~$M$ steps.

One can see that, for $t=0,1,\ldots$,
\begin{multline} \label{eq:estimate1}
\widetilde{G}_1^{(t+1)M+1} - \widetilde{G}_1^{tM+1} =
\biggl(\sum_{i=1}^M \eta_1^{tM+i} -
\mathbf{I} \left(R_1^{tM+1} = 0\right)
\prod_{j=2}^M (1-\alpha_j^{tM+1}) \\
- \sum_{j=2}^M \alpha_1^{tM+j}
\left(\mathbf{I} \left(R_j^{tM+j} =0\right)
\prod_{k=2}^M (1-\alpha_k^{tM+j}) +
\mathbf{I} \left(R_j^{tM+j} > 0\right)
\prod_{k \ge 2, k \neq j} (1-\alpha_k^{tM+j}) \right)\biggr)^+.
\end{multline}
Let~$S_i^{t+1} = \sum_{k=1}^M \xi_i^{tM+k}$ be the total number of arrivals into the red queue of node~$i$ within the consecutive~$M$ time slots, and let

\begin{multline*}
\varphi_{x,1}^{tM+1} = \sum_{i=1}^M \eta_1^{tM+i}
- \mathbf{I} \left(x_1 = 0\right)
\prod_{j=2}^M (1-\alpha_j^{tM+1}) \\
-\sum_{j=2}^M \alpha_1^{tM+j} \left(\mathbf{I}
\left(x_j + S_j^{t+1} = 0\right)
\prod_{k=2}^M (1-\alpha_k^{j}) + \mathbf{I} \left(x_j + S_j > 0\right) \prod_{k \ge 2, k \neq j} (1-\alpha_k^j) \right).
\end{multline*}

Then the RHS of equation \eqref{eq:estimate1} may be estimated from above by random variable $\varphi_{R_1,1}^{tM+1}$.

We prove now that condition~$\widetilde{B1}$ and conditions of Lemma~\ref{suffi} hold for the Markov Chain~$\{R_1^t,...,R_M^t,\widetilde{G}_1^t\}$ in~$M$ steps. Indeed, we have
$$
\mathsf{E} \left(\left|\widetilde{G}_1^{M+1} - \widetilde{G}_1^1\right| |(R_1^1,...,R_M^1,G_1^1) = (r_1,...,r_M,g_1)\right) \le \mathsf{E} \max\left\{\sum_{i=1}^M \eta_1^i, M\right\},
$$
since at most~$\sum_{i=1}^M \eta_1^i$ new messages may arrive in the system and at most~$M$ messages may leave the system. Condition~$\widetilde{B1}$ thus holds. 
Conditions~$(i)$ and~$(ii)$ of Lemma~\ref{suffi} clearly hold for random variables 
$\varphi_{x,1}^{tM+1}$. Take~$C_1 = C_2 = M$, then
 condition~$(iii)$ of Lemma~\ref{suffi} also holds with
$h_i(y_i)=C_1{\mathbf I}(L_{2,i}(y_i)\le C_2)$.
To verify the last condition~$(iv)$, we note
that~$\mathsf P(R^t_{i(t)} = 0) \to 1 - \lambda_R$ as~$t \to \infty$ (this follows from a known general result that, for a stationary Markov chain
$Z_{t+1} =\max (Z_t -1, 0) + \chi_t$ with i.i.d. integer-valued increments $\{\chi_t\}$ such
that ${\mathbf E}\chi_1=c <1$ and ${\mathbf P} (\chi_1 \ge 0 )=1$,
we have with necessity ${\mathbf P} (Z_t=0) = 1-c$, in order to keep our paper self-contained, we provide a proof in Appendix B).
Therefore, we get
\begin{eqnarray*}
\int \mathsf{E} \xi_{x,1}^1 \pi(dx) & = & (1-\lambda_R) p(1-p)^{M-1} \\ & + & \sum_{j=2}^M \left( (1-\lambda_R)p(1-p)^{M-1} + \lambda_R p(1-p)^{M-2} \right) < 0,
\end{eqnarray*}
provided the conditions of the Theorem hold.
\qed

{\sc Proof of instability.}
As was mentioned in the stability proof,~$\mathbf P(R^t_{i(t)} = 0) \to (1 - \lambda_R)$ as $n \to \infty$. We can choose~$t$ so large that~$|\mathbf P(R^t_{i(t)} = 0) - (1 - \lambda_R)| < \delta$ for an arbitrarily small~$\delta > 0$. For simplicity of caculations,
let us assume that~$\mathbf P(R^t_{i(t)} = 0) = (1 - \lambda_R)$ (it will not be difficult for the reader to repeat the same proof with an extra~$\delta$ added and then let~$\delta$ go to~$0$). We prove that, for any $i=1,\ldots, M$,
$G_i^t\to\infty$ with at least a linear speed, i.e.
\begin{equation}\label{instab1}
\liminf_{t\to\infty}G_i^t/t > 0
\quad \mbox{a.s.}
\end{equation}
(note that by the SLLN the speed of growth can not be superlinear). Consider the embedded epochs $M,2M,\ldots, kM,\ldots$ that are the
multiples of $M$. Choose a positive number $N>>1$. Since all
states in the positive $M$-dimensional lattice are communicating,
there exists an a.s. finite (random) time, say, $T \in \{kM, k\ge
0\}$ such that $G_i^T \ge N$ for all $i$. Starting from time $T$,
all coordinates of the process $G^{kM}$ coincide with those of the
auxiliary process $\widetilde{G}^{kM}$ which starts with the same
$\widetilde{G}^T=G^T$ -- until the first time when one of the
coordinates becomes zero. Since, for any $i$, the increments
$\widetilde{G}^{(k+1)M}_{i} - \widetilde{G}^{kM}_i$ form a
stationary ergodic sequence with a positive mean, say $\Delta$ (which is a
difference of the RHS and the LHS in equation~(\ref{eq:cond_aloha})),
$$ \widetilde{G}^{kM}_i/k \to \Delta \quad
\mbox{a.s.} $$ and, for any $\varepsilon >0$, one can choose
$N>>1$ such that $\inf_{l\ge 0} \widetilde{G}^{T+lM} \ge M+1$ with
probability at least $1-\varepsilon$. If one takes $\varepsilon <
1/M$, then all the coordinates of $\widetilde{G}^{kM}$ always stay
above $M$ after time $T$ with probability at least $1-M\varepsilon
>0$. Then the same holds for the coordinates of the process
$G^{kM}$. Since $G^{t+1}_i-G^t_i\ge -1$ a.s., it then follows
that, with probability at least $1-M\varepsilon >0$, the
coordinates of $G^t_i$ stay strictly positive for all $t\ge T$.
Since $\varepsilon $ may be taken as small as possible, (\ref{instab1}) follows.

\qed

Recall that $i(t)= ((t-1) \mod M )+1$.
Introduce a Markov Chain
$$
\left\{\left(\check{R}_1^t,\ldots,
\check{R}_M^t\right), t\ge 0\right\} =
\left\{\left(R_{i(t)}^t, R^t_{i(t+1)},R^t_{i(t+2)},
\ldots, R^t_{i(t+M-1)}\right), t\ge 0\right\}
$$
and another Markov Chain~$(\check{R},\check{G})=
\left\{\left(\check{R}_1^t,\ldots,
\check{R}_M^t,\check{G}_1^t,\ldots,\check{G}_M^t\right),
t\ge 0\right\}$, where a similar interchange of the~$G$-coordinates is also made. In words, the new Markov chain is obtained from the original one with a cyclic change of coordinates such that at every time slot the first coordinate corresponds to the node whose first red message
is scheduled for service (if there is any).
\begin{cor}
Under the assumptions of Theorem~\ref{th:aloha}, the Markov Chain
$(\check{R}, \check{G})$ is Harris ergodic.
\end{cor}

{\sc Proof.} Indeed, the new chain is aperiodic and positive recurrent. The latter follows from the fact that the Lyapunov function used in the proof of Theorem~\ref{th:aloha} is the sum of all coordinates and hence does not change when the order of the coordinates is changed. To infer ergodicity it is then sufficient to show that the state~$(0,...,0)$ is achievable from any other state. To show that, note first that there exists a compact set~$V$ which is positively recurrent for the new chain. This implies that with a positive probability the chain will reach a state~$(r_1,..,r_M,g_1,..g_M) \in V$. Due to compactness of~$V$, there exist finite~$r^0$ and~$g^0$ such that~$\sum_{i=1}^M r_i \le r^0$ and~$\sum_{i=1}^M g_i \le g^0$ for all points from~$V$. It is then clear that there exists a finite time such that with a positive probability all red and all green messages will be transmitted and no new messages will arrive. Hence the state~$(0,...,0)$ may be reached in a finite number of steps from any other state.

\qed

\subsection{Network without a MAC-coordinator} \label{subsec:no_blocking}

Now consider the system which differs from the system described above only
in the following: we assume here that a station which is transmitting a red message can also attempt a
transmission of a green message and, if that happens, these two transmissions collide.

Now the following may happen at a time slot~$t$ regarding red messages:

\begin{itemize}

\item{There is one attempted (and successful) transmission. This
happens when the queue of the red messages of node $i(t)$ is non-empty and there is no attempted
transmission of a green message by the same station (this may happen if either the queue of green
messages is empty, or it is not empty but the station takes a decision not to transmit a green message);}

\item{There is one attempted (and unsuccessful) transmission. This happens when both queues of red and
green messages of station~$i(t)$ are non-empty and the station decides to attempt a transmission of a green message.}

\item{There are no attempted transmissions. This happens when the
queue of the red messages of node $i(t)$ is empty;}

\end{itemize}

The following may occur regarding green messages:

\begin{itemize}

\item{There is one transmission attempt. Then it is successful, if it 
did not originate from node~$i(t)$
or if it did originate from node~$i(t)$ but its queue of red messages is empty;}

\item{No transmission attempted;}

\item{Two or more transmissions attempted. In this case all these
transmissions fail due to the collision.}

\end{itemize}

We are still assuming that the ALOHA algorithm is used to govern the behaviour of the green queues
of all stations. For this system, the following holds.

\begin{theorem} \label{th:aloha2}
The system is stable if $\lambda_R < 1-p$ and
\begin{equation} \label{eq:cond_aloha2}
\begin{cases}
\lambda_G < \left(1-\frac{\lambda_R}{1-p}\right) p, \quad \text{if} \quad M=1, \cr
\lambda_G < \frac{\lambda_R}{1-p} (M-1) p (1-p)^{M-1} + \left(1-\frac{\lambda_R}{1-p}\right)
Mp(1-p)^{M-1}, \quad \text{if} \quad M>1.
\end{cases}
\end{equation}
\end{theorem}

{\sc Proof} of Theorem \ref{th:aloha2}
may be given following the lines of the proof of Theorem~\ref{th:aloha}. The only difference is
that one needs to bound the Markov Chain representing the state of the system under consideration by a Markov Chain representing the state of the system where each station with an empty green queue
transmits a ``dummy" message with probability~$p$. 

Note that in this case the first component
(counting the number of red messages) of the initial Markov Chain describing the state of the system is not
a Markov Chain itself, so Theorem~\ref{th:stability_multi} can not be applied directly. However,
the Markov Chain used to bound the initial one from above has the first component which is a Markov
 Chain itself, and the use of Theorem~\ref{th:stability_multi} is justified.

{\bf Acknowledgement.} The authors would like to thank Maxim Grankin for his simulation work that numerically confirmed the correctness of stability results in the case of a network without a MAC-coordinator. 

\section*{Appendix A}

{\sc Proof of \eqref{c_uni1}-\eqref{c_uni2}}.
Assume $p<1$ (where~$p$ is from the minorisation condition, see Assumption~$A3$)-- otherwise the result is obvious.
Consider, for simplicity, the case $m=1$ only (the general
case requires an extra technical work which is not essential).
For $x\in V$, consider the standard splitting identity
$$
{\mathbf P}_x (X^1 \in \cdot ) =
p\mu (\cdot ) + (1-p)
\frac{{\mathbf P}_x (X^1 \in \cdot ) - p\mu(\cdot )}{1-p}
$$
and denote the fraction in the RHS by ${\mathbf Q}_x(\cdot )$
which is a probability measure.
We know that sigma-algebra ${\mathcal B_X}$ is countably generated.
Therefore, one can define two measurable functions
$f,g: {\mathcal X} \times [0,1] \to {\mathcal X}$ such that if $U$ is a r.v.
uniformly distributed in $[0,1]$, then
$f(x,U)$ has distribution ${\mathbf P}_x(\cdot )$, for $x\in {\mathcal X}$,
and $g(x,U)$ has
distribution ${\mathbf Q}_x(\cdot )$, for $x\in V$ -- see, e.g.,
\cite{Kifer} for background.

Introduce three sequences
of mutually independent r.v.'s, each of which is i.i.d.:\\
1) a sequence of $0-1$-valued r.v.'s $\beta_t$, with
common distribution
${\mathbf P} (\beta_t=1)=p=1-{\mathbf P} (\beta_t=0)$,\\
2) a sequence $\{U_t\}$ of uniformly distributed in $[0,1]$
r.v.'s, and\\
3) a sequence $\{W_t\}$ of ${\mathcal X}$-valued r.v.'s with
common distribution $\mu$.

Then Markov chain $X^t$ may be represented as a stochastic recursive
sequence (SRS):
\begin{equation}\label{SRS}
X^{t+1} = \left( \beta_t W_t + (1-\beta_t) g(X^t,U_t)\right)
{\mathbf I} (X^t\in V) + f(X^t,U_t) {\mathbf I}(X^t \in {\mathcal X}
\setminus V).
\end{equation}
The pairs $(X^t,\beta_t)$ also form a time-homogeneous
Markov chain. Start with $X^0=x\in V$. Let $T_0=0$ and, for $k\ge 1$,
$$
T_{k} = \min \{t> T_{k-1} : \ X^t \in V\}.
$$
Further, let
$$
\theta = \min \{ k\ge 0: \  \beta_{T_k}=1\}.
$$
Then $\theta$ has a geometric distribution with parameter $p$,
${\mathbf P} (\theta = k)= p(1-p)^k$, $k=0,1,\ldots$. Let
$\kappa = T_{\theta}+1$. Note that
r.v. $X^{\kappa}$ has distribution $\mu$ and that, for $x\in V$,
$$
{\mathbf E}_x \kappa \le s_0 {\mathbf E} \theta + 1 =: C
$$
where  $s_0$ is from \eqref{eq:def_recurrence_2}.
Clearly, $C$ does not depend on $x\in V$. Then, in particular,
random variables $\kappa$ are uniformly bounded in probability:
$$
\sup_{x\in A} {\mathbf P}_x (\kappa >t) \to 0, \ t\to\infty,
$$
by Markov inequality.

Let now $\{\bar{\beta}_t\},
\{\bar{U}_t\},$ and $\{\bar{W}_t\}$ be three other i.i.d. sequences
which do not depend on the first
three sequences.
Consider a stationary Markov chain~$\{\bar{X}^t\}$ which is defined
as follows: $\bar{X}^0$ has distribution $\pi$ and does not depend
on all r.v.'s defined earlier, and
\begin{equation}\label{SRS2}
\bar{X}^{t+1} = \left( \bar{\beta}_t \bar{W}_t +
(1-\bar{\beta}_t)
g(\bar{X}^t,\bar{U}_t)\right)
{\mathbf I} (\bar{X}^t\in V) +
f(\bar{X}^t,\bar{U}_t)
{\mathbf I}(\bar{X}^t \in {\mathcal X}
\setminus V).
\end{equation}
Due to independence of the two SRS's, r.v. $\bar{X}^{\kappa}$ has
distribution $\pi$.
Finally, let
$$
\gamma = \min \{ j\ge 0 : \
X^{\kappa+j} \in A, \bar{X}^{\kappa
+j}\in A, \beta_{\kappa+j}=1 \}.
$$
By aperiodicity, $\gamma$ is finite a.s. Also, it does not
depend on $\kappa$ and, therefore,
its distribution does not depend on $x$.
Then one can define another sequence $\widehat{X}^n$
by
$$
\widehat{X}^t = \bar{X}^t {\bf I} (t\le \kappa + \gamma )
+X^t {\bf I} (t>\kappa + \gamma )
$$
and find that, first, since~$\kappa+\gamma$ is a stopping time, $\{\bar{X}^t\}$ and~$\{\widehat{X}^t\}$ have the same distribution and, in particular, $\{\widehat{X}^t\}$ is also a stationary
Markov chain (see, e.g. \cite[p.34]{Thorisson} or \cite{Lindvall}), second, random variables $\nu = \kappa + \gamma$ are uniformly bounded
in probability, i.e. \eqref{c_uni2} holds, and, third,
r.v.'s $\nu$ satisfy \eqref{c_uni1}.
\qed

{\bf Remark.} The first intention of the authors was to find this result
in \cite{Thorisson}. However, Hermann Thorisson has confirmed
that it is not there, but he is thinking to have it in the second edition
of the book (with a complete proof, for any $m\ge 1$).

\section*{Appendix B}

\begin{lemma}
Consider a stationary Markov chain $Z_{t+1} =\max (Z_t -1, 0) + \chi_t$ with a stationary and ergodic sequence of integer-valued increments $\{\chi_t\}$ such
that ${\mathbf E}\chi_1=c <1$ and ${\mathbf P} (\chi_1 \ge 0 )=1$. Then ${\mathbf P} (Z_t=0) = 1-c$.
\end{lemma}

{\sc Proof.}

Let us rewrite the evolution equation as
$$
Z_{t+1} = Z_t - 1 + {\bf I}\{Z_t=0\} + \chi_t,
$$
where~${\bf I}\{\cdotp\}$ is the indicator function. With a positive integer~$T$ we now have
$$
Z_{t+T} = Z_t - T + \sum_{i=0}^{T-1} {\bf I}\{Z_{t+i}=0\} + \sum_{i=0}^{T-1} \chi_{t+i}.
$$
If we now divide both sides of the latter expression by~$T$, take expectations and let~$T$ tend to infinity, we will get
$$
\lim_{T \to \infty} \frac{{\mathbf E} Z_{t+T}}{T} = -1 + {\mathbf P}(Z_t = 0) + c, 
$$
where we used the Strong Law of Large Numbers and the stationarity of the Markov chain~$\{Z_t\}$. Using stationarity again, we conclude that the LHS of the latter equality is equal to~$0$, and the proof is complete. \qed

\end{document}